\documentclass[12pt,reqno]{article}

\usepackage{graphicx}
\usepackage{epstopdf}
\usepackage{amsmath}
\usepackage{amsthm}
\usepackage{latexsym,bm}
\usepackage{amssymb}
\usepackage{indentfirst}
\newtheorem{thm}{Theorem}[section]
\newtheorem{cor}[thm]{Corollary}
\newtheorem{lem}[thm]{Lemma}

\numberwithin{equation}{section}
\usepackage{graphicx}

\begin{document}

\title{{\bf
On the complexity of isometric immersions of hyperbolic spaces in
any  codimension
%$\mathbb H^m \to \mathbb R^n$
%of  hyperbolic
%$m$-space into  $\mathbb R^n$
%There is no simple way of embedding isometrically hyperbolic
%$\text{m}$-space in  $\mathbb R^n$
%The complexity of the isometric immersions of  hyperbolic
%$m$-space into  $\mathbb R^n$
}}

\author{{\bf F. Fontenele\thanks{Work partially supported by CNPq
(Brazil)}\, and F. Xavier}}

\date{}
\maketitle

\begin{quote}
\small {\bf Abstract}. Although the Nash theorem solves the
isometric embedding problem, matters are inherently more involved
if one is further seeking an embedding that is well-behaved from the
standpoint of submanifold geometry. More generally, consider a
Lipschitz map $F:M^m\to \mathbb R^n$, where $M^m$ is a Hadamard
manifold whose curvature lies between negative constants. The main
result of this paper is that $F$ must perform a substantial
compression: For every $r>0$ and integer $k\geq 2$ there exist $k$
geodesic balls of radius $r$ in $M^m$ that are arbitrarily far
from each other, but whose images under $F$ are bunched together
arbitrarily close in the Hausdorff sense of $\mathbb R^n$. In
particular, every isometric embedding $\mathbb H^m\to \mathbb R^n$
of hyperbolic space must have a complex  asymptotic behavior,
regardless of how high the codimension is. Hence, there is no
truly simple way to realize $\mathbb H^m$ isometrically inside any
Euclidean space.

\end{quote}

\section{Introduction.}

The Nash embedding theorem (\cite {K}, \cite{N}), to the effect
that any Riemannian manifold $(M^m, g)$ can be isometrically
embedded as a bounded subset of some high dimensional Euclidean
space, represents a landmark in Riemannian geometry. Questions
about the smoothness of the isometric immersion, connections with
topology and partial differential equations, as well as the
smallest dimension of the receiving space, have also attracted
considerable attention over the years (\cite{G}, \cite {GR},
\cite{HH}).

Exploring the surrounding landscape further, one comes across the
natural idea of establishing the existence of an isometric
immersion that, from the perspective of global submanifold geometry, is
as well-behaved as the intrinsic geometry allows.

For instance, if $(M^m,g)$ is non-compact but \it complete\rm, one
might aim for the existence of a \it proper \rm  isometric
immersion (or embedding) $F:M^m\to \mathbb R^n$, for some $n>m$. A
more refined problem, albeit vaguely stated, would be to produce a
proper isometric embedding whose behavior at infinity is as tame
as possible.

A suitable testing ground for these ideas is $\mathbb H^m$, the
complete simply-connected space of constant sectional curvature
$-1$.  For reasons that were not clear,
it  has been   rather difficult to realize the Nash theorem in the case of $\mathbb H^m$,
namely to produce explicit
isometric embeddings of $\mathbb H^m$ into \it some \rm Euclidean space.

One such result is due to  Henke and Nettekoven \cite{HN},  where $\mathbb H^m$ is  \it properly \rm isometrically
embedded in $\mathbb R^{6m-6}$ as a smooth complete graph over an
$m$-dimensional subspace (see also \cite{B},
\cite{He}). How can the existence of such a proper
isometric embedding be reconciled with the fact that hyperbolic
space $\mathbb H^m$ is much larger at infinity than any Euclidean
space? What is the role of the codimension?

The main finding of this paper is that, in order to accommodate
the different orders of growth at infinity, any  isometric
embedding $F:\mathbb H^m \to \mathbb R^n$ - regardless of
regularity, dimension or codimension - must exhibit a high degree
of asymptotic complexity, which is expressed in a precise
quantitative fashion.

In fact, as we shall see below,  our arguments can be implemented in
the broader context of Lipschitz maps $M^m \to \mathbb R^n$, where $M^m$ is  a   Hadamard
manifold with curvature bounded away from zero.

Before stating the main result, we explain its geometric meaning
in an informal way:

\vskip10pt

\noindent ($\dagger$) \it For every $r>0$ and integer $k\geq 2$
there exist  $k$ geodesic balls of radius $r$ in $\mathbb H^m$ that are arbitrarily
far apart, but whose images  under  $F$ are
arbitrarily Hausdorff-close in $\mathbb R^n$. \rm

\vskip10pt
When $F$ is the lift to $\mathbb H^m$ of an
isometric immersion of a compact hyperbolic manifold and $r$ is
large enough, the images of the balls actually coincide with
$F(\mathbb H^m)$.  If $F$ is a proper isometric embedding, as in
\cite{HN}, and $r$, $k$ are arbitrary but fixed, one obtains from $(\dagger$) the
following ``dynamical"  picture, that helps in the visualization of  the embedding:
\vskip10pt
\noindent ($\dagger \dagger$)
\it There is a sequence of
configurations of $k$ balls of radius $r$ in $\mathbb H^m$, the
distance between any two balls in a configuration going to
infinity, such that the isometric images  under  $F$ of the $k$ balls form a sequence of ``stacks"  in
$\mathbb R^n$, each one with $k$ ``layers", with the  property that  the stacks tend to
infinity in $\mathbb R^n$  while their thickness tends  to
zero.\rm
\vskip10pt

Hence,  despite the fact that
$\mathbb H^m$ is a simple space, any isometric embedding $\mathbb
H^m \hookrightarrow \mathbb R^n$, proper or not, must be rather
complex, as the submanifold in $\mathbb R^n$ fails to stabilize at
infinity.

The  phenomenon described above helps to explain  why    isometric embeddings of hyperbolic spaces
into Euclidean spaces  are  so hard to produce explicitly, since  any candidate  for an isometric embedding must  exhibit  \it a  priori \rm  a specific complex asymptotic behavior.

There are    simple models of hyperbolic geometry
that retain   some features   of Euclidean geometry, for instance the  classical models
of   Poincar\'e, Lobachevsky,
Minkowski-Lorentz and Cayley-Klein.
On the other hand, our results  reveal   that,
despite Nash's theorem,  there is  no  truly
simple   way to  actually realize hyperbolic spaces
isometrically \rm  inside any Euclidean space,  no matter  how high
the codimension is allowed to be.

We observe that matters are much simpler when the roles of the
spaces are reversed. For instance, horospheres provide
well-behaved examples of isometric embeddings $\mathbb
R^n\hookrightarrow \mathbb H^m, \;m=n+1$.

As mentioned before, our main result holds for maps that are  more general
than isometric immersions $\mathbb H^m\to\mathbb R^n$. Its formal
statement reads as follows: \vskip10pt
\begin{thm}\label{Teorema 1}
Let $m\geq 2$, $n\geq 1$ be integers, $M^m$ a Hadamard manifold whose
curvature is bounded above by a negative constant, and $F:M^m\to
\mathbb R^n$ a Lipschitz map. Then, for every $r>0$, $\epsilon \in
(0,1)$ and integer $k\geq 2$, there are points $p_1,\dots,p_k\in
M^m$ for which the geodesic balls $B(p_i,r)$ satisfy, for all
distinct $i,j\in\{1,\dots,k\}:$
\vskip5pt \noindent i) The
Riemannian distance between $B(p_i, r)$ and $B(p_j, r)$ is at
least $\epsilon^{-1}$. \vskip3pt \noindent ii) The Euclidean
distance between $F(B(p_i,r))$ and $F(B(p_j,r))$ is at most
$\epsilon$.
\vskip5pt \noindent If the curvature of $M^m$ lies
between negative constants,
then i) and iii) below hold:
\vskip3pt \noindent iii)  The Hausdorff distance between
$F(B(p_i,r))$ and $F(B(p_j,r))$ is at most $\epsilon$.\end{thm}

\vskip10pt

The proof of Theorem \ref{Teorema 1}, to be given in the next
section, is based on a careful study of the interplay between the
asymptotic growth of some specially defined packings of geodesic
balls in the strongly curved Hadamard manifold $M^m$,
and the massive compression they must undergo under the action of a Lipschitz
map that takes values in some Euclidean space.

The conceptual remarks below are meant to shed some light on
Theorem \ref{Teorema 1}, vis-a-vis the nature of the known
examples of isometric immersions $\mathbb H^m\to\mathbb R^n$. Call
a map $G:M\to N$ between non-compact complete Riemannian manifolds
\it strongly proper \rm if, for every pair of sequences $(x_n)$,
$(y_n)$ in $M$, the following condition is fulfilled:
$$d(x_n,y_n)\to\infty\implies\liminf
d(G(x_n),G(y_n))>0.$$

As the terminology indicates, a strongly proper map can easily be
seen to be proper in the usual sense. On the other hand, the
following example shows that not every proper immersion is
strongly proper. Consider $A\subset \mathbb R$ given by the
disjoint union, over all integers $k\geq 1$, of the open intervals
$(k-1/(k+1), k+1/(k+1))$, and let $f:\mathbb R\to\mathbb R$ be a
smooth function satisfying $\lim_{k\to \infty} f(k)=\infty$,
$f(x)>0$ for $x\in A$, $f(x)=0$ for $x\in\mathbb R-A$. Define
$g:\mathbb R^2\to\mathbb R$ by $g(x,y)=f(x)$. It is easy to see
that the (complete, properly embedded) graph $S\subset \mathbb
R^3$ of $g$ is not strongly proper. Indeed, the points
$p_k=(k-1/(k+1),0,0)$ and $q_k=(k+1/(k+1),0,0)$ lie in $S$ and
satisfy $||p_k-q_k||=2/(k+1)\to 0$. But, on the other hand, one
can argue using the symmetries of $S$ that $d_S(p_k,q_k)\geq
2f(k)\to\infty$.

The following is a direct consequence of Theorem \ref{Teorema 1}:
\vskip5pt
\begin{cor}\label{nice}
Let $M^m$ be a Hadamard manifold whose curvature is bounded above
by a negative constant, and $F:M^m\to\mathbb R^n$ a Lipschitz map.
Then, for every integer $k\geq 2$ there are $k$ sequences
$(x_l^{(i)})$ in $M^m$, $1\leq i\leq k$, such that, for all
distinct $i,j\in\{1,\dots,k\}$, one has
$$\lim_{l\to\infty}d(x_l^{(i)},x_l^{(j)})\to\infty \;\;\; \text{and}\;\;\;  \lim_{l\to\infty}||F(x_l^{(i)})-F(x_l^{(j)})||= 0.$$
In particular, $F$ is not strongly proper ($k=2$).
\end{cor}

Some proper isometric immersions $\mathbb R^m\to\mathbb R^n$
between flat Euclidean spaces, for instance the totally geodesic
ones, are strongly proper. But observe that the graph $S$ of the
example above is flat, proper, but not strongly proper.  Isometric immersions of hyperbolic spaces into $\mathbb
R^n$, on the other hand,  behave quite differently. Although \cite{HN} provides
examples of  proper  isometric embeddings $\mathbb H^m
\hookrightarrow \mathbb R^{6m-6}$, Corollary \ref{nice} implies:

\vskip5pt
\begin{cor} \label{nice too}
There are examples of proper  isometric immersions $F:\mathbb
H^m\to\mathbb R^n$, but no such $F$ can be strongly proper.
\end{cor}

Needless to say, the full force of Theorem \ref{Teorema 1}
provides much more information about isometric immersions $\mathbb
H^m \to \mathbb R^n$ than Corollary \ref{nice too}.

For the sake of completeness, we  mention  the well-known problem  that $\mathbb H^m$ cannot be $C^2$ isometrically immersed
in $\mathbb R^{2m-1}$, although this conjecture is not the focus
of the present work (indeed, our results are valid in arbitrary codimension).  For background on this problem, as well as
related works, see \cite{E}, \cite{F}, \cite{H}, \cite{Mo},
\cite{Ni} - \cite {X}).

We would like to stress that since Theorem
\ref{Teorema 1} holds for Lipschitz functions, it  is conceivable  that it
may  be of use in other problems in geometric analysis, besides  isometric immersions.

Given the somewhat general nature of our arguments, we expect
Theorem \ref{Teorema 1} to admit  formulations in other settings as well,
provided that  there is  a notion of hyperbolicity that can be played
against the idea of  polynomial    growth.

In conclusion, despite the fact that the Nash theorem solves the isometric
embedding problem, matters are inherently more involved from the
standpoint of the geometry of submanifolds. Indeed, as it will be seen in
this paper, in some cases there are global obstructions at work
that preclude the existence of isometric embeddings with tame
asymptotic behavior,  \it regardless of the codimension.\rm

The authors would like to thank F. Ledrappier for useful
conversations regarding some aspects of this work.

\section{Lipschitz maps and asymptotic densities.}

\vskip10pt

This section contains the proof of Theorem \ref{Teorema 1},
presented after some preparatory material. Given a complete
non-compact $m$-dimensional Riemannian manifold $M^m$, $R\in
(0,\infty)$, $C\in (0,R)$, and $p\in M^m$, denote by
$\#(p,C,R;M^m)$ the maximum number of disjoint metric balls of
radius $C$ that are contained in the open ball $B(p,R)$.

\vskip5pt
\begin{lem}\label{cresc}
Let $M^m$ be a Hadamard manifold with
curvature bounded away from zero. Then, for all $p\in M^m$ and
$C>0$, $\#(p,C,R;M^m)$ grows exponentially as $R\to\infty$.
\end{lem}

In the special case $M^m=\mathbb H^m$, a non-computational proof
of this result can be given using the fact that there are compact
hyperbolic manifolds $P^m$ with an arbitrarily large injectivity
radius, together with the exponential growth of the fundamental
group of $P^m$. Since we were unable to locate a reference for
Lemma \ref{cresc} in the case of variable curvature, a detailed proof  will be provided.

The following standard result (\cite{dC}, \cite{CE}) will be used
in Lemmas \ref{cresc} and \ref {EpsilonNet}: \vskip5pt
\begin{lem}\label{CorRauch}
Let $M^n$ and $\widetilde M^n$ be Riemannian manifolds and suppose
that $\widetilde K_{\widetilde x}(\widetilde\sigma)\geq
K_x(\sigma)$, for all $x\in M$, $\widetilde x\in\widetilde M$,
$\sigma\subset T_xM$, $\widetilde\sigma\subset T_{\widetilde
x}\widetilde M$. Let $p\in M$, $\widetilde p\in\widetilde M$ and
fix a linear isometry $\iota:T_pM\to T_{\widetilde p}\widetilde
M$. Let $r>0$ such that $\textnormal{exp}_p|_{B_r(0)}$ is a
diffeomorphism and $\textnormal{exp}_{\widetilde p}|_{\widetilde
B_r(0)}$ is non-singular. Let $c:[0,a]\to
\textnormal{exp}_{p}\big(B_r(0)\big)\subset M$ be smooth and
define  $\widetilde c:[0,a]\to\textnormal{exp}_{\widetilde
p}\big(\widetilde B_r(0)\big)\subset\widetilde M$ by $ \widetilde
c(s)=\textnormal{exp}_{\widetilde p}\circ
\iota\circ\textnormal{exp}_p^{-1}\big(c(s)\big),\;s\in [0,a].$
Then $l(c)\geq l(\widetilde c)$.
\end{lem}
In order to prove Lemma \ref {cresc} we may assume, without loss
of generality, that the sectional curvature $K$ of $M^m$ is at
most $-1$. As before, denote by $\mathbb H^m$ the $m$-dimensional
hyperbolic space and fix $p\in M^m$, $\tilde p\in\mathbb H^m$.
Since $\text{exp}_p:T_pM\to M$, $\text{exp}_{\tilde p}:T_{\tilde
p}\,\mathbb H^m\to\mathbb H^m$ are diffeomorphisms, for any fixed
linear isometry $\iota:T_pM\to T_{\tilde p}\,\mathbb H^m$ the map
$\phi:=\text{exp}_{\tilde
p}\circ\iota\circ\text{exp}_p^{-1}:M^m\to\mathbb H^m$ is also a
diffeomorphism. Moreover, by Lemma \ref{CorRauch},
\begin{eqnarray}\label{Hadam1A}
d_M(x,y)\geq d_{\mathbb H^m}\big(\phi(x),\phi(y)\big),\;\;\;x,y\in
M^m.
\end{eqnarray}

Our strategy will be to identify a suitable two-dimensional
surface with the property that the maximum number of disjoint
balls of radius $C$ that are contained in $B(p,R)$, and whose
centers lie in the said surface, already grows exponentially as
$R \to \infty$.

Given $R>2C$, let $\alpha\in(0,\pi/2)$ be such that
\begin{eqnarray}\label{Hadam1}
\sin\alpha=\frac{\sinh C}{\sinh (R-C)}\,\cdot
\end{eqnarray}
Let $k$ be largest positive integer that satisfies $k\alpha\leq
\pi-\alpha$, so that
\begin{eqnarray}\label{Hadam2}
k>\frac{\pi-\alpha}{\alpha}-1.
\end{eqnarray}
Let $u,w\in T_pM$ be orthogonal unit vectors and, for each integer
$0\leq j\leq k$, set
$$
v_j=\cos(2j\alpha)u+\sin(2j\alpha)w.
$$

We claim that, for all $0\leq i<j\leq k$, the smallest angle
$\angle(v_i,v_j)$ between $v_i$ and $v_j$ is at least $2\alpha$.
In fact, taking the inner product of $v_i$ and $v_j$ we obtain
$$
\cos\angle(v_i,v_j)=\cos(2j\alpha-2i\alpha)=\cos\big(2\pi-(2j\alpha-2i\alpha)\big).
$$
If $2j\alpha-2i\alpha\leq\pi$, then
$\angle(v_i,v_j)=(j-i)2\alpha\geq 2\alpha$. On the other hand, if
$2j\alpha-2i\alpha>\pi$ we have
$$
\angle(v_i,v_j)=2\pi-(j-i)2\alpha\geq 2\pi-2k\alpha.
$$
That $\angle(v_i,v_j)\geq 2\alpha$ is valid also in this case is
an immediate consequence of the above inequality and our choice of
$k$.

For $0\leq j\leq k$, consider the geodesic in $M$ given by
$\gamma_j(t)=\text{exp}_p(tv_j)$,  and let $p_j=\gamma_j(R-C)$.

\vskip10pt

\begin{lem} \label{Claim} $d_M(p_i,p_j)\geq 2C$ for all distinct
$i,j\;$ in $\{0,1,...,k\}$. \end{lem}

\noindent To prove  Lemma \ref{Claim}, for each $j$ such that $0\leq j\leq k$ consider
the geodesic in $\mathbb H^m$ defined by
$\tilde\gamma_j(t)=\text{exp}_{\tilde p}\big(t\,\iota(v_j)\big)$,
and let $\tilde p_j=\tilde\gamma_j(R-C)$. Since $\phi(p_j)=\tilde
p_j$, it follows from (\ref{Hadam1A}) that
\begin{eqnarray}\label{Hadam3}
d_M(p_i,p_j)\geq d_{\mathbb H^m}(\tilde p_i,\tilde
p_j),\;\;\;i,j\in\{0,1,...,k\},\;i\neq j.
\end{eqnarray}

We can assume that $\angle(v_i,v_j)$ is strictly less than $\pi$, otherwise
$\tilde\gamma_i$ and $\tilde\gamma_j$ would be opposite geodesics,
and so
$$
d_{\mathbb H^m}(\tilde p_i,\tilde p_j)=2(R-C)>2C.
$$

At this point in the proof we need to invoke a classical formula
in hyperbolic trigonometry, the so-called hyperbolic law of sines
\cite[p. 432]{Ma}. This formula states that if a triangle in the
hyperbolic plane has sides of lengths $a$, $b$, $c$, and the
corresponding opposite angles have measures $\lambda$, $\mu$,
$\nu$, then
$$
\frac{\sinh a}{\sin\lambda}=\frac{\sinh b}{\sin\mu}=\frac{\sinh
c}{\sin\nu}.
$$

Let $\tilde q$ be the midpoint of the segment $\tilde p_i\tilde
p_j$, so that the triangle $\tilde p\tilde q\tilde p_i$ has a
right angle at $\tilde q$ and an angle of measure $\frac{1}{2}
\angle(v_i,v_j)$ at $\tilde p$. Applying the hyperbolic law of
sines to the triangle $\tilde p\tilde q\tilde p_i$, one obtains
\begin{eqnarray}\label{Hadam6}
\sinh\big(d_{\mathbb H^m}(\tilde q,\tilde
p_i)\big)=\sinh(R-C)\sin\left(\frac{\angle(v_i,v_j)}{2}\right).
\end{eqnarray}
Since, by the previous claim, $\angle(v_i,v_j)\geq 2\alpha$, it
follows from (\ref{Hadam1}) and (\ref{Hadam6}) that
\begin{eqnarray}\label{Hadam6a}
\sinh\big(d_{\mathbb H^m}(\tilde q,\tilde
p_i)\big)\geq\sinh(R-C)\sin\alpha=\sinh C.
\end{eqnarray}
Lemma \ref{Claim} now follows from  (\ref{Hadam3}) and (\ref{Hadam6a}):
\begin{eqnarray}\label{Hadam6b}
d_M(p_i,p_j)\geq d_{\mathbb H^m}(\tilde p_i,\tilde
p_j)=2d_{\mathbb H^m}(\tilde q,\tilde p_i)\geq 2C.
\end{eqnarray}

We now resume the proof of Lemma \ref{cresc}. As $d_M(p,p_j)=R-C$,
it follows from the triangle inequality that $B(p_j,C)\subset
B(p,R)$ for $0\leq j\leq k$. Moreover, by Lemma \ref{Claim},
$B(p_i,C)\cap B(p_j,C)=\emptyset$ for all distinct
$i,j\in\{0,1,...,k\}$. From (\ref{Hadam1}) and (\ref{Hadam2}), we
then obtain
\begin{eqnarray}\label{Hadam8}
\#(p,C,R;M^m)&\geq&
k+1>\frac{\pi-\alpha}{\alpha}\nonumber\\&=&\frac{\sin\alpha}{\alpha}\,\frac{\pi-\alpha}{\sin\alpha}
\nonumber\\&=&\frac{\sin\alpha}{\alpha}\,\frac{\pi-\alpha}{\sinh
C}\,\sinh
(R-C)\nonumber\\&>&\frac{1}{2}\,\frac{\sin\alpha}{\alpha}\,\frac{\pi}{\sinh
C}\,\sinh (R-C).
\end{eqnarray}
Since, by (\ref{Hadam1}), $\alpha\to 0$ as $R\to\infty$, it
follows from (\ref{Hadam8}) that $\#(p,C,R;M^m)$ grows
exponentially with $R$, for every $C>0$ fixed. This concludes the
proof of Lemma \ref{cresc}.\qed

\vskip10pt

\begin{lem}\label{EpsilonNet}
Let $M^m$ be a Hadamard manifold whose sectional curvature $K$ is
bounded from below. Then, for all $\rho>0$ and $\delta>0$ there
exist a positive integer $l=l(\rho,\delta)$ and maps
$\sigma_1,...,\sigma_l:M\to M$ such that, for all $p\in M$,

\vskip5pt

\noindent i) $\sigma_j(p)\in B(p,\rho),\;\;1\leq j\leq l$,

\vskip5pt

\noindent ii)
$B(p,\rho)\subset\bigcup_{j=1}^lB(\sigma_j(p),\delta)$.
\end{lem}

To prove the lemma, set $b=\inf_MK>-\infty$ and let $\widetilde M$
be the complete simply-connected $m$-dimensional Riemannian
manifold with constant sectional curvature $b$. Fix $\widetilde
q\in\widetilde M$ and let $\widetilde q_1,...,\widetilde q_l\in
B(\widetilde q,\rho)$ be such that
\begin{eqnarray}\label{EpsilonNet1}
B(\widetilde q,\rho)\subset\bigcup_{j=1}^lB(\widetilde
q_j,\delta).
\end{eqnarray}
For $q\in M$ fixed, consider orthonormal bases $\{\widetilde
v_1,...,\widetilde v_m\}$ and $\{v_1,...,v_m\}$ of $T_{\widetilde
q}\widetilde M$  and $T_qM$, respectively. For each $p\in M$,
$p\neq q$, let $\{V_1(p),...,V_m(p)\}$ be the (orthonormal) basis
of $T_pM$ obtained by the parallel transport of $v_1,...,v_m$
along the (unique) geodesic joining $q$ to $p$. Consider also the
linear isometry $\iota_p:T_{\widetilde q}\widetilde M\to T_pM$
satisfying $\iota_p(\widetilde v_i)=V_i(p)$, and define a
diffeomorphism $\phi_p:\widetilde M\to M$ by
$\phi_p=\text{exp}_p\circ\iota_p\circ\text{exp}_{\widetilde
q}^{-1}.$

For all $p\in M$ and $j\in\{1,...,l\}$,  set
$\sigma_j(p)=\phi_p(\widetilde q_j)\in M$. Since
$\phi_p(B(\widetilde q,\rho))=B(p,\rho)$, we have $\sigma_j(p)\in
B(p,\rho)$ whenever $1\leq j\leq l$. Given $x\in B(p,\rho)$, we
obtain from (\ref{EpsilonNet1}) that $\phi_p^{-1}(x)\in
B(\widetilde q_j,\delta)$ for some $j$, $1\leq j\leq l$. Applying
Lemma \ref{CorRauch} with the roles of $M$ and $\widetilde M$
interchanged, one has
\begin{eqnarray*}
d_{\widetilde M}\big(\phi_p^{-1}(x),\phi_p^{-1}(y)\big) \geq d_M(x,y),
\end{eqnarray*}
and so $ d_M(x,\sigma_j(p))\leq d_{\widetilde
M}(\phi_p^{-1}(x),\widetilde q_j)<\delta, $ which establishes
ii).\qed

\vskip20pt

With these preliminaries out of the way, we are ready to begin the
proof of Theorem \ref{Teorema 1}. To this end, fix $p_0\in M^m$
and let $C$ be a positive number to be specified later. For each
$R>C$, take a collection $\widehat{\mathcal C}_{R,C}$ of disjoint
balls of radius $C$ inside the ball $B(p_0,R)\subset M^m$ such
that $|\widehat{\mathcal C}_{R,C}|=\#(p_0,C,R;M^m)$.

According to Lemma \ref{cresc}, the cardinality
$|\widehat{\mathcal C}_{R,C}|$ grows exponentially
as $R\to\infty$. In particular, one has
\begin{eqnarray}\label{outro limsupA}
\lim_{R\to\infty}\frac{|\widehat{\mathcal
C}_{R,C}|}{R^n}=\infty.
\end{eqnarray}

Consider, for each $R>C$, a subfamily $\mathcal C_{R,C}$ of
$\widehat{\mathcal C}_{R,C}$ satisfying: \vskip10pt \noindent a)
If $|\mathcal C_{R,C}|>1$ and $p,q$ are centers of distinct balls
in $\mathcal C_{R,C}$, then $||F(p)-F(q)||\geq\frac{1}{C}.$
\vskip5pt \noindent b) $\mathcal C_{R,C}$ is maximal with respect
to property a).

\vskip10pt

Writing $D(q,t)$ for the Euclidean ball in $\mathbb R^n$ of radius
$t$ and center $q$, we observe that if $|\mathcal C_{R,C}|>1$
and $B(q_i,C)$, $B(q_j,C)$ are distinct balls in $\mathcal
C_{R,C}$, then
\begin{eqnarray}\label{interA}
D\left(F(q_i),\frac{1}{3C}\right)\cap D\left(F(q_j),
\frac{1}{3C}\right)=\emptyset.
\end{eqnarray}
Indeed, if (\ref{interA}) were to fail, the triangle inequality
would imply $||F(q_i)-F(q_j)||<\frac{2}{3C}$, contradicting a)
above.

Denote by $L$ the Lipschitz constant of $F$. Since
\begin{eqnarray*}
||F(q_i)-F(p_0)||\leq Ld(q_i,p_0)<LR
\end{eqnarray*}
we have, for all $x\in D\big(F(q_i),1/3C\big)$,
\begin{eqnarray*}
||x-F(p_0)||\leq ||x-F(q_i)||+||F(q_i)-F(p_0)||<\frac{1}{3C}+LR.
\end{eqnarray*}
As a consequence,
\begin{eqnarray}\label{inclusionA}
\bigcup_{q_i} D\left(F(q_i),\frac{1}{3C}\right)\subset
D\left(F(p_0),LR+\frac{1}{3C}\right),
\end{eqnarray}
where $q_i$ runs over the centers of all balls in $\mathcal
C_{R,C}$.

An individual ball $D(q,t)$ in $\mathbb R^n$ has volume $c(n)t^n$,
the explicit value of the constant $c(n)$ being unimportant for
our current purposes. The volume of each ball
$D\big(F(q_i),1/3C\big)$ is a fixed constant, say $\lambda_0$. There
are $|\mathcal C_{R,C}|$ such balls in $\mathbb R^n$ and, as
observed in (\ref{interA}), they are pairwise disjoint. Thus, the
volume of the union in (\ref{inclusionA}) is $\lambda_0 |\mathcal
C_{R,C}|$ and, furthermore,
\begin{eqnarray}\label{compA}
\lambda_0 |\mathcal C_{R,C}| \leq
c(n)\big(LR+1/3C\big)^n.\nonumber
\end{eqnarray}
In particular,
\begin{eqnarray}\label{limsupA}
\limsup_{R\to\infty}\frac{|\mathcal C_{R,C}|}{R^n}<\infty.
\end{eqnarray}

It follows from (\ref{outro limsupA}) and (\ref{limsupA}) that,
for all sufficiently large $R$, say $R>R_0$, the inclusion
$\mathcal C_{R,C}\subset\widehat {\mathcal C}_{R,C}$ is proper.

Consider any ball $B(p,C)$ from  $\widehat{\mathcal
C}_{R,C}-\mathcal C_{R,C}$. The family $\{B(p,C)\}\cup\mathcal
C_{R,C}$ consists of disjoint balls of radius $C$ and so, by a)
and the maximality of $\mathcal C_{R,C}$ that was stipulated in
b), one can select a (not necessarily unique) ball $B(q,C)$ in
$\mathcal C_{R,C}$ such that
\begin{eqnarray}\label{menor queA}
||F(p)-F(q)||<\frac{1}{C}.
\end{eqnarray}
Any such assignment $B(p,C)\leadsto B(q,C)$ gives rise to a map
\begin{eqnarray}
\Theta_{R,C}: \widehat {\mathcal C}_{R,C}-\mathcal C_{R,C}\to
\mathcal C_{R,C},\;\;R>R_0.\nonumber
\end{eqnarray}

We claim that when $R$ tends to infinity, the cardinality of \it
some \rm  fiber of $\Theta_{R,C}$ becomes larger than any
specified integer $j$.

Indeed, if not,
\begin{eqnarray*}
|\widehat{\mathcal C}_{R,C}|\nonumber&=&|\mathcal
C_{R,C}|+|\widehat{\mathcal C}_{R,C}-\mathcal
C_{R,C}|\\\nonumber
&\leq&|\mathcal C_{R,C}|+j|\Theta_{R,C}(\widehat{\mathcal C}_{R,C}-\mathcal C_{R,C})|\\
&\leq&(1+j)|\mathcal C_{R,C}|,
\end{eqnarray*}
contradicting (\ref{outro limsupA}) and (\ref{limsupA}).

Hence, by (\ref{menor queA}) and the previous assertion about the
fibers of $\Theta_{R,C}$, there are points $q,p_1,\dots,p_k\in M^m$ such that
\begin{eqnarray*}
\min _{1\leq i,j\leq k,\;i\neq j}d(p_i,p_j)\geq 2C\;\;\;
\text{and}\;\;\;\max_{1\leq i\leq
k}||F(p_i)-F(q)||<\frac{1}{C}\,\cdot
\end{eqnarray*}
In particular,
\begin{eqnarray*}
\max_{1\leq i,j\leq k}||F(p_i)-F(p_j)||<\frac{2}{C}\,\cdot
\end{eqnarray*}
Choosing $C=2(2r\epsilon+1)/\epsilon$, we obtain (ii) in the statement of Theorem \ref{Teorema 1}. To see that
(i) is also valid with this choice of $C$, observe that, for all
$x\in B(p_i,r)$ and $y\in B(p_j,r)$,
\begin{eqnarray*}
d(x,y)\geq d(p_i,p_j)-2r\geq
2C-2r>\frac{2r\epsilon+1}{\epsilon}-2r=\frac{1}{\epsilon}\,.
\end{eqnarray*}

We now  move on  to the second half of the theorem, and assume that the
sectional curvature of $M^m$ is bounded from above and below by
negative constants.

By Lemma \ref{EpsilonNet}, with $\rho=r$ and $\delta=\frac{\epsilon}{2L}$, there exist a positive integer $l$ and
maps $\sigma_1,...,\sigma_l:M\to M$ such that, for all $p\in M$,
\begin{eqnarray}\label{coverAA}
\sigma_j(p)\in B(p,r),\;\;1\leq j\leq l,
\end{eqnarray}
\begin{eqnarray}\label{coverA}
B(p,r)\subset\bigcup_{j=1}^lB\big(\sigma_j(p),\epsilon/2L\big).
\end{eqnarray}

In order to control the Hausdorff distance, we introduce the following
augmentation of the map $F$:
\begin{eqnarray}\label{chapeu}
\widehat F:M^m\to\mathbb R^n\times \cdots \times \mathbb R^n=\mathbb R^{nl},\;\;\;\widehat
F(p)=\big(F(\sigma_1(p)),...,F(\sigma_l(p))\big).
\end{eqnarray}

Fix $p_0\in M^m$,  and let $C$ be a positive number to be specified
later. As before, for each $R>C$ denote  by $\widehat{\mathcal
C}_{R,C}$ a collection of disjoint balls of radius $C$ that are contained in   the  ball
$B(p_0,R)\subset M^m$  and satisfy  $|\widehat{\mathcal C}_{R,C}|=\#(p_0,C,R;M^m)$.

By Lemma
\ref{cresc}, one has
\begin{eqnarray}\label{outro limsup}
\lim_{R\to\infty}\frac{|\widehat{\mathcal
C}_{R,C}|}{R^{nl}}=\infty.
\end{eqnarray}
Consider, for each $R>C$, a subfamily $\widetilde{\mathcal
C}_{R,C}$ of $\widehat{\mathcal C}_{R,C}$ satisfying: \vskip5pt
\noindent a) If $|\widetilde{\mathcal C}_{R,C}|>1$ and $p,q$ are
centers of distinct balls in $\widetilde{\mathcal C}_{R,C}$, then
$||\widehat F(p)-\widehat F(q)||\geq\frac{1}{C}.$ \vskip5pt
\noindent b) $\widetilde{\mathcal C}_{R,C}$ is maximal with
respect to property a).

\vskip10pt

From this point on, we write $D(q,t)$ for the ball in $\mathbb
R^{nl}$ of radius $t$ and center $q$. If $|\widetilde{\mathcal
C}_{R,C}|>1$ and $B(q_i,C)$, $B(q_j,C)$ are distinct balls in
$\widetilde{\mathcal C}_{R,C}$ then, by a) above,
\begin{eqnarray}\label{inter}
D\left(\widehat F(q_i),\frac{1}{3C}\right)\cap D\left(\widehat
F(q_j), \frac{1}{3C}\right)=\emptyset.
\end{eqnarray}
From
\begin{eqnarray*}
||F(\sigma_s(p))-F(\sigma_s(p_0))||\leq
Ld(\sigma_s(p),\sigma_s(p_0))\leq L\big
[d(\sigma_s(p),p)+d(p,p_0)+d(p_0,\sigma_s(p_0)\big]
\end{eqnarray*}
and (\ref{coverAA}), one obtains
\begin{eqnarray*}
||F(\sigma_s(p))-F(\sigma_s(p_0))||<L\big(d(p,p_0)+2r\big),
\end{eqnarray*}
\begin{eqnarray*}
||\widehat F(p)-\widehat F(p_0)||^2=\sum_{s=1}^l
||F(\sigma_s(p))-F(\sigma_s(p_0))||^2<
lL^2\big(d(p,p_0)+2r\big)^2.
\end{eqnarray*}
Then, for all $x\in D(\widehat F(q_i),1/3C))$ we have
\begin{eqnarray}
||x-\widehat F(p_0)||\leq ||x-\widehat F(q_i)||+||\widehat
F(q_i)-\widehat F(p_0)||<\frac{1}{3C}+\sqrt{l}L(R+2r).\nonumber
\end{eqnarray}
As a consequence,
\begin{eqnarray}\label{inclusion}
\bigcup_{q_i} D\left(\widehat F(q_i),\frac{1}{3C}\right)\subset
D\left(\widehat F(p_0),\sqrt{l}L(R+2r)+\frac{1}{3C}\right),
\end{eqnarray}
where $q_i$ runs over the centers of all balls in $\widetilde{\mathcal C}_{R,C}$.

Denoting by $\overline{\lambda}$ the volume of each ball $D(\widehat
F(q_i),1/3C)$, it follows from (\ref{inter}) and (\ref{inclusion})
that
\begin{eqnarray}\label{comp}
\overline{\lambda} |\widetilde{\mathcal C}_{R,C}|\leq
c(nl)\big(\sqrt{l}L(R+2r)+1/3C\big)^{nl},\nonumber
\end{eqnarray}
where $c(nl)$ is the volume of the unit ball in $\mathbb R^{nl}$.
In particular,
\begin{eqnarray}\label{limsup}
\limsup_{R\to\infty}\frac{|\widetilde{\mathcal C}_{R,C}|}{R^{nl}}<\infty.
\end{eqnarray}

It follows from (\ref{outro limsup}) and (\ref{limsup}) that, for
all sufficiently large $R$, say $R>R_0$, the inclusion
$\widetilde{\mathcal C}_{R,C}\subset\widehat{\mathcal C}_{R,C}$ is
proper. Hence, by a) and the maximality of $\widetilde{\mathcal
C}_{R,C}$ that was stipulated in b), for any ball $B(p,C)$ from
$\widehat{\mathcal C}_{R,C}-\widetilde{\mathcal C}_{R,C}$ one can
select a (not necessarily unique) ball $B(q,C)$ in
$\widetilde{\mathcal C}_{R,C}$ such that
\begin{eqnarray}\label{menor queAa}
||\widehat F(p)-\widehat F(q)||<\frac{1}{C}\,.
\end{eqnarray}
Any such assignment $B(p,C)\leadsto B(q,C)$ gives rise to a map
\begin{eqnarray*}
\widetilde{\Theta}_{R,C}:
\widehat {\mathcal C}_{R,C}-\widetilde{\mathcal C}_{R,C}\to\widetilde{\mathcal
C}_{R,C},\;\;R>R_0.
\end{eqnarray*}
If all fibers of  $\widetilde{\Theta}_{R,C}$ had size at most a fixed integer $j$, an  argument  similar to the discussion following (\ref{menor queA}) would show that
$|\widehat {\mathcal C}_{R,C}| \leq (1+j)  \widetilde{\mathcal C}_{R,C}$,  contradicting  (\ref{outro limsup}) and (\ref{limsup}).

Hence, by (\ref{menor queAa}) there
are points $q,p_1,\dots,p_k\in M^m$ for which
\begin{eqnarray*}
\min _{1\leq i,j\leq k,\;i\neq j}d(p_i,p_j)\geq 2C\;\;\;
\text{and}\;\;\;\max_{1\leq i\leq k}||\widehat F(p_i)-\widehat
F(q)||<\frac{1}{C}\,\cdot
\end{eqnarray*}
In particular,
\begin{eqnarray}\label{newA}
\max_{1\leq i,j\leq k}||\widehat F(p_i)-\widehat
F(p_j)||<\frac{2}{C}\,\cdot
\end{eqnarray}

For $x\in B(p_i,r)$ and $y\in B(p_j,r)$, one has
\begin{eqnarray}
d(x,y)\geq d(p_i,p_j)-2r\geq 2C-2r,\nonumber
\end{eqnarray}
and so, under the hypothesis that the curvature lies between
negative constants, Theorem \ref{Teorema 1} i) follows by choosing
any $C$ satisfying
\begin{eqnarray*}
C\geq\frac{2r\epsilon+1}{2\epsilon}.
\end{eqnarray*}

Given $x\in B(p_i,r)$,  (\ref{coverA}) implies that there
is $s\in\{1,\dots,l\}$ with $d(x,\sigma_s(p_i))<\epsilon/2L$,
so that, by the Lipschitz condition,
\begin{eqnarray}\label{new2}
||F(x)-F(\sigma_s(p_i))||<\frac{\epsilon}{2}\,.\nonumber
\end{eqnarray}
Since, by (\ref{chapeu}) and (\ref{newA}),
\begin{eqnarray}
||F(\sigma_s(p_i))-F(\sigma_s(p_j))||\leq ||\widehat
F(p_i)-\widehat F(p_j)||<\frac{2}{C},\nonumber
\end{eqnarray}
we have
\begin{eqnarray}\label{intermediate}
||F(x)-F(\sigma_s(p_j))||<\frac{\epsilon}{2}+\frac{2}{C}.
\end{eqnarray}
Finally, choosing
\begin{eqnarray*}
C\geq\max\left\{\frac{2r\epsilon+1}{2\epsilon},\frac{4}{\epsilon}\right\}
\end{eqnarray*}
one sees from (\ref{intermediate}) that $F(x)$ lies in the
$\epsilon$-neighborhood of the set $F(B(p_j,r))$. As $x\in
B(p_i,r)$ is arbitrary, $F(B(p_i,r))$ is contained  in the
$\epsilon$-neighborhood of $F(B(p_j,r))$. Theorem \ref{Teorema 1}
iii) now follows by reversing the roles of $i$ and $j$ in the
argument above.

%\begin{flushleft}
$$
\begin{array}{lcccccccccl}
\text{Francisco Fontenele}            &&&&&&&&& & \text{Frederico Xavier}\\
\text{Departamento de Geometria}      &&&&&&&&& & \text{Department of Mathematics}\\
\text{Universidade Federal Fluminense}&&&&&&&&& & \text{University of Notre Dame}\\
\text{Niter\'oi, RJ, Brazil}          &&&&&&&&& & \text{Notre Dame, IN, USA}\\
\text{fontenele@mat.uff.br}           &&&&&&&&& & \text{fxavier@nd.edu}\\
\end{array}
$$
%\end{flushleft}
\end{document}